\documentclass[10pt]{article}

\usepackage{a4wide}
\usepackage{amssymb}
\usepackage{amsfonts}
\usepackage{amsmath}
\input xy
\xyoption{arrow} \xyoption{matrix}

\date{}

\newtheorem{proposition}{Proposition}[section]
\newtheorem{theorem}[proposition]{Theorem}
\newtheorem{lemma}[proposition]{Lemma}

\newtheorem{corollary}[proposition]{Corollary}

\def\der{\partial }

\def\nFM0{{\nu }_{F,M_0}}
\def\nFN0{{\nu }_{F,N_0}}
\def\nGN0{{\nu }_{G,N_0}}

\def\N0{ {\bf N}_0 }

\def\t{\otimes}

\def\ra{\rightarrow}

\def\Xpm{X^{\pm }}

\def\s{\sigma}

\def\l1{{\lambda}_1}

\def\a{\alpha}
\def\a0{ {\alpha }_0}
\def\a1{ {\alpha }_1}

\def\l{\lambda}

\def\nFGM0{{\nu }_{F,G,M_0}}

\def\nFN0{{\nu}_{F,N_0}}

\def\sm{{\sigma}^m}

\def\sm1{{\sigma}^{-1}}

\def\smtp1{{\sigma}^{-t+1}}

\def\S1{S^{-1}}

\def\Xpm1{X^{\pm 1}_1}

\def\sPM1{{\sigma }^{\pm 1}}
\def\sMP1{{\sigma }^{\mp 1 }}

\def\di{{\rm d.ind}}

\def\L{\Lambda}

\def\Ytm1{Y^{t-1}}
\def\Yim1{Y^{i-1}}

\def\CM{{\cal M}}
\def\CN{{\cal N}}

\def\ass{{\rm ass}}

\def\bQ{\overline{Q}}

\def\ker{ {\rm ker } }

\def\SL2Z{ {\rm SL}_2({\bf Z}) }

\def\Gp1{ G^{1 , 1 } }
\def\P11{ P^{-1 , 1 } }
\def\Pp1{ P^{1 , 1 } }

\def\nCLsr{{}^\nu\kern-2pt {\cal L}^{\sigma , \rho  }}
\def\nP{{}^\nu \kern-2pt P}
\def\nL{{}^\nu\kern-2pt L}
\def\nLL{{}^\nu\kern-2pt \Lambda}
\def\nPsr{{}^\nu\kern-2pt P^{\sigma , \rho  }}
\def\nLsr{{}^\nu\kern-2pt L^{\sigma , \rho  }}
\def\nuCL{{}^\nu\kern-2pt  {\cal L}}
\def\nCLsr{{}^\nu\kern-2pt {\cal L}^{\sigma , \rho  }}
\def\nCL1m{{}^\nu\kern-2pt {\cal L}^{-1 , 1  }}
\def\x1nu{x^\frac{1}{\nu}}
\def\xm1nu{x^{-\frac{1}{\nu}}}

\def\rad{{\rm rad}}

\def\CN{{\cal N}}
\def\ra{\rightarrow }

\def\CB{{\cal B}}

\def\CT{{\cal T}}

\def\CC{ {\cal C}}

\def\nAM0{{\nu }_{{\cal A},M_0}}
\def\nAN0{{\nu }_{{\cal A},N_0}}

\def\bI{\overline{I}}
\def\bp{\overline{p}}
\def\bR{\overline{R}}

\def\bQ{\overline{Q}}

\def\bq{\overline{q}}

\def\ga{\mathfrak{a}}
\def\gb{\mathfrak{b}}

\def\gn{\mathfrak{n}}
\def\gm{\mathfrak{m}}
\def\gp{\mathfrak{p}}
\def\gq{\mathfrak{q}}

\def\SL{{\rm SL}}

\def\tga{\widetilde{\mathfrak{a}}}

\def\di!{\frac{\der^i}{i!}}
\def\dik!{\frac{\der^k_i}{k!}}

\def\N{\mathbb{N}}

\def\0{\overline{0}}
\def\1{\overline{1}}

\def\Ln1{\L_{n,\overline{1}}}

\def\a1{a_{\overline{1}}}

\def\bs{\overline{s}}

\def\S{\Sigma}

\def\vn1{\overrightarrow{n-1}}

\def\Min{{\rm Min}}

\def\mJ{\mathbb{J}}
\def\mI{\mathbb{I}}

\def\ann{{\rm ann}}

\def\K1{{\rm K}_1}

\def\hmI1{\widehat{\mI_1}}
\def\tmI1{\widetilde{\mI_1}}
\def\tmJ1{\widetilde{\mJ_1}}
\def\hB1{\widehat{B_1}}
\def\hCB1{\widehat{\CB_1}}

\def\bS{\overline{S}}

\def\Den{{\rm Den}}

\def\Ore{{\rm Ore}}

\def\Den{{\rm Den}}

\def\Ass{{\rm Ass}}

\def\maxDen{{\rm max.Den}}

\def\llrad{{\rm l.lrad}}

\def\OCC{\overline{\CC}}
\def\br{\overline{r}}
\def\bc{\overline{c}}

\def\bs{\overline{s}}
\def\bt{\overline{t}}
\def\ga{\mathfrak{a}}

\def\tor{{\rm tor}}

\def\gll{\mathfrak{l}}

\def\tCC{\widetilde{\mathcal{C}}}
\def\dCC{\mathcal{C}^\dagger}
\def\tQ{\widetilde{Q}}
\def\ts{\widetilde{\sigma}}
\def\grR{{\rm gr}\, R}
\def\grQ{{\rm gr}\, Q}

\begin{document}

\author{V. V. \  Bavula 
}

\title{Criteria for a ring to have a left Noetherian  left quotient ring}

\maketitle

\begin{abstract}
 Two criteria are given for a ring to have a left Noetherian left quotient ring (this was an open problem since 70's). It is proved that each such ring has only {\em finitely many} maximal left denominator sets.

$\noindent $

 {\em Key Words:  Goldie's Theorem, 
 the left quotient ring  of a ring, the largest left quotient ring of a ring, a maximal left denominator set, 
 the prime radical.}

 {\em Mathematics subject classification
 2010: 15P50,  16P60,  16P20, 16U20.}

$${\bf Contents}$$
\begin{enumerate}
\item Introduction.
\item  Properties of rings with a left Noetherian left quotient ring.
 \item Proof of Theorem \ref{6Jan13}.
and  the Second Criterion (via the associated graded ring).
 \item Finiteness of $\maxDen_l(R)$ when $Q(R)$ is a left Noetherian ring.
\end{enumerate}
\end{abstract}


\section{Introduction}

In this paper, module means a left module, and the following notation is fixed:
\begin{itemize}
\item  $R$ is a ring with 1;

\item   $\CC = \CC_R$  is the set of {\em regular} elements of the ring $R$ (i.e. $\CC$ is the set of non-zero-divisors of the ring $R$);
\item   $Q=Q(R):=Q_{l,cl}(R):= \CC^{-1}R$ is the {\em left quotient ring}  (the {\em classical left ring of fractions}) of the ring $R$ (if it exists) and $Q^*$ is the group of units of $Q$;
\item   $\gn =\gn_R$ is the  prime radical of $R$,  $\nu\in \N \cup \{ \infty \}$ is its {\em nilpotency degree} ($\gn^\nu \neq 0$ but $\gn^{\nu +1}=0$) and $\CN_i := \gn^i/ \gn^{i+1}$ for $i\in \N$;
\item   $\bR := R/ \gn$ and $\pi: R\ra \bR$, $r\mapsto \br =r+\gn$;
\item   $\OCC := \CC_{\bR}$ is the  set of regular elements of the ring $\bR$ and $\bQ := \OCC^{-1}\bR$ is its left quotient ring;
  \item $\tCC :=\pi (\CC )$, $\tQ :=\tCC^{-1}\bR$ and $\dCC := \CC_{\tQ }$ is the set of regular elements of the ring $\tQ$;


\item $S_l=S_l(R)$ is the {\em largest left Ore} of $R$ that consists of regular elements
and
$Q_l=Q_l(R):=S_l(R)^{-1}R$ is the {\em largest left quotient ring} of $R$  {\cite[Theorem 2.1]{larglquot}}; \item $\Ore_l(R):=\{ S\, | \, S$ is a left Ore set in $R\}$; \item
$\Den_l(R):=\{ S\, | \, S$ is a left denominator set in $R\}$.

\end{itemize}
{\bf Criteria for  a ring to have a left Noetherian  left quotient ring}.
The aim of the paper is to give two criteria for a ring $R$ to have a left Noetherian  left  quotient ring (Theorem \ref{6Jan13} and Theorem \ref{27Jan13}).
 The case when $R$ is a {\em semiprime} ring  is a very easy special case.

 \begin{theorem}\label{4Oct14}
Let $R$ be a semiprime ring. Then the following statements are equivalent.
\begin{enumerate}
\item $Q(R)$ is a left Noetherian ring.
\item $Q(R)$ is a semisimple ring.
\item $R$ is a semiprime  left Goldie ring.
\item $Q_l(R)$ is a left Noetherian ring.
\item $Q_l(R)$ is a semisimple ring.
\end{enumerate}
If one of the equivalent conditions holds then $S_l(R) = \CC_R$ and $Q= Q_l(R)$. In particular, if the left quotient ring $Q$ (resp. $Q_l(R)$)  is not a semisimple ring then the ring $Q$ (resp. $Q_l(R)$) is not left Noetherian.
\end{theorem}
{\em Example}, \cite{Bav-intdifline}. The ring $\mI_1:= K\langle x, \frac{d}{dx}, \int \rangle$ of polynomial integro-differential operators over a field $K$ of characteristic zero  is a semiprime ring but not left Goldie (as it contains infinite direct sums of non-zero left ideals). Therefore, the largest left quotient ring $Q_l(\mI_1)$ is not a left Noetherian ring (moreover, the left quotient ring $Q(\mI_1)$ does not exists). The ring $Q_l(\mI_1)$ and $S_l(\mI_1)$ were  described explicitly in \cite{Bav-intdifline}.

The first criterion for a ring to have a left Noetherian left quotient ring is below, its proof is given in Section \ref{PTHA6Jan13}. 

\begin{theorem}\label{6Jan13}
 Let $R$ be a ring. The following statements are equivalent.

\begin{enumerate}
\item  The left quotient ring $Q(R)$ of $R$  is a left Noetherian ring.
\item
\begin{enumerate}
\item $\tCC \subseteq \OCC$.
\item $\tCC\in \Ore_l(\bR )$.
\item $\tQ = \tCC^{-1}\bR$ is a left Noetherian ring.
\item $\gn$ is a nilpotent ideal  of the ring $R$.
\item The $\tQ$-modules $\tCC^{-1}\CN_i$, $i=1, \ldots , \nu $, are   finitely generated  (where $\nu $ is the nilpotency degree of $\gn$ and  $\CN_i:=\gn^i/\gn^{i+1}$).
    \item For each element $\bc\in \tCC$, the left $\bR$-module $\CN_i/ \CN_i\bc$ is $\tCC$-torsion for $i=1, \ldots , \nu $.
       \end{enumerate}
\end{enumerate}
\end{theorem}

{\it Remark}. The conditions (a) and (b) above imply that the set $\tCC$ is a left denominator set of the ring $R$, and so the ring $\tQ$ exists.

The powers of the prime radical $\gn$, $\{ \gn^i\}_{i\geq 0}$, form a descending filtration on $R$. Let $\grR = \bR\oplus\gn / \gn^2\oplus\cdots$ be the associated graded ring. The second criterion is given in terms of  properties of the ring $\grR$, its proof is given in Section \ref{PTHA6Jan13}. 

\begin{theorem}\label{27Jan13}
Let $R$ be a ring. The following statements are equivalent.
\begin{enumerate}
\item The ring $R$ has a left Noetherian left quotient ring $Q$.
\item The set $\tCC$ is a left denominator set of the ring $\grR$, $\tCC \subseteq \OCC$, $\tCC^{-1}\grR $ is a left Noetherian ring and $\gn$ is  a nilpotent ideal.
\end{enumerate}
If one of the equivalent conditions holds then $\grQ\simeq \tCC^{-1} \grR$ where $\grQ :=\tQ\oplus\gn_Q / \gn_Q^2\oplus\cdots$ is the associated graded ring with respect to the prime radical filtration. In particular, the ring $\grQ$ is a left Noetherian ring.
\end{theorem}

{\bf Finiteness of the set of maximal left denominators for a ring with a left Noetherian left quotient ring}.
 For a ring $R$, the set $\maxDen_l(R)$ of maximal left denominator sets  (with respect to $\subseteq$) is a {\em non-empty} set, \cite{larglquot}. It was proved that the set $\maxDen_l(R)$ is a finite set if the left quotient ring $Q(R)$ of $R$ is a semisimple ring, \cite{Bav-Crit-S-Simp-lQuot},  or a left Artinian ring,  \cite{Bav-LocArtRing}. The next theorem extends this result for a larger class of rings that includes the class of left Noetherian ring for which the left quotient ring exists. 
\begin{theorem}\label{aA5Oct14}
Let $R$ be a ring such that its left quotient ring  $ Q(R)$ is a left Noetherian ring. Then $|\maxDen_l(R)|<\infty$. Moreover, $|\maxDen_l(R)|\leq s=|\maxDen_l(\bR )|$ where $\bQ\simeq  \prod_{i=1}^s\bQ_i$ and $\bQ_i$ are simple Artinian rings (see Theorem \ref{A6Jan13}.(3)).
\end{theorem}

The next corollary follows at once from Theorem \ref{aA5Oct14}.

\begin{corollary}\label{AA5Oct14}
If $R$ is a left Noetherian ring such that its left quotient ring $Q(R)$ exists (i.e. $\CC\in \Ore_l(R)$)  then $|\maxDen_l(R)|<\infty$. 
\end{corollary}

The next corollary is an explicit description of the set $\maxDen_l(R)$ for a ring $R$  with a left Noetherian left quotient ring $Q(R)$. 

\begin{corollary}\label{xaA5Oct14}
Let $R$ be a ring such that its left quotient ring $Q(R)$ is a left Noetherian ring. For each $i=1,\ldots , s$, let $p_i: R\ra \bQ_i$ be the natural projection (see (\ref{piRRQ}) and Theorem \ref{aA5Oct14}), $\bQ_i^*$ be the group of units of the simple Artinian ring $\bQ_i$, $S_i'$ be the largest element (w.r.t. $\subseteq$) of the set $D_i=\{ S'\in \Den_l(R)\, | \, p_i(S')\subseteq \bQ_i^*\}$. Then 
\begin{enumerate}
\item $\maxDen_l(R)$ is the set of maximal elements  (w.r.t. $\subseteq$) of the set $\{ S_1', \ldots , S_s'\}$.  
\item For all $i=1,\ldots , s$, $\CC\subseteq S_i'$.
\item The rings $S_i'^{-1}R$ are left Noetherian where $i=1, \ldots , s$.  
\end{enumerate}
\end{corollary}


{\bf A criterion for a left Noetherian ring to have a left Noetherian left quotient ring}. In \cite{Goldie-1969AspRTh}, Goldie posed a problem of deciding when a Noetherian ring possesses a Noetherian quotient ring. In \cite{Stafford-1982NoethQuotRing},  Stafford obtained a criterion to determine when a Noetherian ring is its own quotient ring. In \cite{Chat-Hajarn-1994}, Chatters and Hajarnavis obtained  necessary and sufficient conditions for a Noetherian ring which is a finite module over its centre to have a quotient ring.

 The next corollary  follows from Theorem \ref{6Jan13}, Theorem \ref{27Jan13}  and the fact that the prime radical of a left Noetherian ring is a nilpotent ideal.
\begin{corollary}\label{a20Nov14}
Let $R$ be a left Noetherian ring. The following statements are equivalent.
\begin{enumerate}
\item The set $\CC =\CC_R$ is a left Ore set (i.e. the left quotient ring $Q(R)=\CC^{-1}R$ of $R$  exists). \item $\tCC \in \Den_l(\bR , 0)$ and for  each element $\bc\in \tCC$ the left $\bR$-module $\CN_i/ \CN_i\bc$ is $\tCC$-torsion for $i=1, \ldots , \nu $.
\item $\tCC\subseteq \CC$ and $\tCC\in \Den (\grR)$.
\end{enumerate}
\end{corollary}

{\it Proof}. $(1\Leftrightarrow 2)$ and   $(1\Leftrightarrow 3)$ follow from Theorem \ref{6Jan13} and Theorem \ref{27Jan13}, respectively. $\Box $

\begin{itemize}
\item {\em Corollary \ref{a8Jan13} is a criterion for a ring to have a left Noetherian left quotient  ring $Q$ such that the factor ring $Q/\gn_Q$ is a semisimple ring (or $\tQ$ is a semisimple ring; or $\tCC = \OCC$; or $\CC = \pi^{-1}(\OCC )$).}
 
 \item {\em Theorem \ref{A23Nov14} is a criterion for a ring $R$ to have a left Noetherian ring such that $|\maxDen_l(R)|=|\maxDen_l(\bR )|$ (recall  that, in general,  $|\maxDen_l(R)|\leq |\maxDen_l(\bR )|$, Theorem \ref{aA5Oct14}).}
 \end{itemize}

{\bf The set of maximal denominator sets for a commutative ring}.  The next theorem describes the set of maximal denominator sets (i.e. maximal Ore sets) for a commutative ring. 

\begin{theorem}\label{23Nov14}
Let $R$ be a commutative ring and $\Min (R)$ be the set of minimal prime ideals of the ring $R$.  Then $\maxDen (R) = \{ S_\gp := R\backslash \gp \, | \, \gp \in \Min (R)\}$. 
\end{theorem}

The paper is organized as follows. In Section \ref{PRLNLQ}, many properties of a ring $R$ with a left Noetherian left quotient ring are proven (Theorem \ref{A6Jan13}). In particular, the implication $(1\Rightarrow 2)$ of Theorem \ref{A6Jan13} is proven.

In Section \ref{PTHA6Jan13},  proofs of Theorem \ref{6Jan13} and Theorem \ref{27Jan13} are given. In Section \ref{FMDQRLN}, Theorem \ref{aA5Oct14} is proven.


\section{Properties of rings with a left Noetherian left quotient ring}\label{PRLNLQ}

In this section, we establish many properties of rings with left Noetherian left quotient ring (Theorem \ref{A6Jan13}). In particular, the implication $(1\Rightarrow 2)$ of Theorem \ref{6Jan13} is proven (Theorem \ref{A6Jan13}.(2)). A criterion is given (Corollary \ref{a8Jan13}) for the ring $\tQ$ to be a semisimple ring or a left Artinian ring.


At the beginning  of the section, we collect necessary results that are used in the proofs of this paper. More results on localizations of rings (and some of the missed standard definitions) the reader can find in \cite{Jategaonkar-LocNRings}, \cite{Stenstrom-RingQuot} and \cite{MR}.  In this paper the following notation will remain fixed:



\begin{itemize}
\item $S_\ga=S_\ga (R)=S_{l,\ga }(R)$
 is the {\em largest element} of the poset $(\Den_l(R, \ga ),
\subseteq )$ and $Q_\ga (R):=Q_{l,\ga }(R):=S_\ga^{-1} R$ is  the
{\em largest left quotient ring associated with} $\ga$, $S_\ga $
exists (Theorem 
 2.1, \cite{larglquot});
\item In particular, $S_0=S_0(R)=S_{l,0}(R)$ is the largest
element of the poset $(\Den_l(R, 0), \subseteq )$ and
$Q_l(R):=S_0^{-1}R$ is the {\em largest left quotient ring} of $R$, \cite{larglquot};

\end{itemize}

{\bf The largest regular left Ore set and the largest left
quotient ring of a ring}. Let $R$ be a ring. A {\em
multiplicatively closed subset} $S$ of $R$ or a {\em
 multiplicative subset} of $R$ (i.e. a multiplicative sub-semigroup of $(R,
\cdot )$ such that $1\in S$ and $0\not\in S$) is said to be a {\em
left Ore set} if it satisfies the {\em left Ore condition}: for
each $r\in R$ and
 $s\in S$, $ Sr\bigcap Rs\neq \emptyset $.
Let $\Ore_l(R)$ be the set of all left Ore sets of $R$.
  For  $S\in \Ore_l(R)$, $\ass (S) :=\{ r\in
R\, | \, sr=0 \;\; {\rm for\;  some}\;\; s\in S\}$  is an ideal of
the ring $R$.

$\noindent $

A left Ore set $S$ is called a {\em left denominator set} of the
ring $R$ if $rs=0$ for some elements $ r\in R$ and $s\in S$ implies
$tr=0$ for some element $t\in S$, i.e. $r\in \ass (S)$. Let
$\Den_l(R)$ be the set of all left denominator sets of $R$. For
$S\in \Den_l(R)$, let $S^{-1}R=\{ s^{-1}r\, | \, s\in S, r\in R\}$
be the {\em left localization} of the ring $R$ at $S$ (the {\em
left quotient ring} of $R$ at $S$). In Ore's method of localization one can localize {\em precisely} at left denominator sets.

In general, the set $\CC$ of regular elements of a ring $R$ is
neither left nor right Ore set of the ring $R$ and as a
 result neither left nor right classical  quotient ring ($Q_{l,cl}(R):=\CC^{-1}R$ and
 $Q_{r,cl}(R):=R\CC^{-1}$) exists.
 There  exists the largest (w.r.t. $\subseteq$) 
 regular left Ore set $S_0= S_{l,0} = S_{l,0}(R)$, \cite{larglquot}. This means that the set $S_{l,0}(R)$ is an Ore set of
 the ring $R$ that consists
 of regular elements (i.e., $S_{l,0}(R)\subseteq \CC$) and contains all the left Ore sets in $R$ that consist of
 regular elements. Also, there exists the largest regular  right (resp. left and right) Ore set  $S_{r,0}(R)$ (resp. $S_{l,r,0}(R)$) of the ring $R$.
 In general, all the sets $\CC$, $S_{l,0}(R)$, $S_{r,0}(R)$ and $S_{l,r,0}(R)$ are distinct, for example,
 when $R= \mI_1= K\langle x, \der , \int\rangle$  is the ring of polynomial integro-differential operators  over a field $K$ of characteristic zero,  \cite{Bav-intdifline}. In  \cite{Bav-intdifline},  these four sets are found for $R=\mI_1$.

$\noindent $

{\it Definition}, \cite{Bav-intdifline}, \cite{larglquot}.    The ring
$$Q_l(R):= S_{l,0}(R)^{-1}R$$ (respectively, $Q_r(R):=RS_{r,0}(R)^{-1}$ and
$Q(R):= S_{l,r,0}(R)^{-1}R\simeq RS_{l,r,0}(R)^{-1}$) is  called
the {\em largest left} (respectively, {\em right and two-sided})
{\em quotient ring} of the ring $R$.

$\noindent $

 In general, the rings $Q_l(R)$, $Q_r(R)$ and $Q(R)$
are not isomorphic, for example, when $R= \mI_1$, \cite{Bav-intdifline}. 

 Small and Stafford \cite{Small-Stafford-1982} have shown that any (left and right) Noetherian ring $R$ possesses a uniquely determined set of prime ideals $P_1, \ldots , P_n$ such that $\CC_R=\cap_{i=1}^n \CC (P_i)$, an irreducible intersection, where $\CC (P_i) :=\{ r\in R\, |\, r+P_i\in \CC_{R/P_i}\}$. Michler and M\"{u}ller \cite{Michler-Muller-1984MaxOre} mentioned that the ring $R$ contains a unique maximal (left and right) Ore set of regular elements
 $S_{l,r,0}(R)$ and called the ring $Q(R)$ the {\em total quotient ring} of $R$. For certain Noetherian rings, they described the set $S_{l,r,0}(R)$ and  the ring $Q(R)$. For the class of affine Noetherian PI-rings, further generalizations were given by M\"{u}ller in \cite{Muller-1985-AffPI}.

 The next
theorem gives various properties of the ring $Q_l(R)$. In
particular, it describes its group of units.


\begin{theorem}\label{4Jul10}
(\cite{larglquot}.)
\begin{enumerate}
\item $ S_0 (Q_l(R))= Q_l(R)^*$ {\em and} $S_0(Q_l(R))\cap R=
S_0(R)$.
 \item $Q_l(R)^*= \langle S_0(R), S_0(R)^{-1}\rangle$, {\em i.e. the
 group of units of the ring $Q_l(R)$ is generated by the sets
 $S_0(R)$ and} $S_0(R)^{-1}:= \{ s^{-1} \, | \, s\in S_0(R)\}$.
 \item $Q_l(R)^* = \{ s^{-1}t\, | \, s,t\in S_0(R)\}$.
 \item $Q_l(Q_l(R))=Q_l(R)$.
\end{enumerate}
\end{theorem}

{\bf The maximal denominator sets and the maximal left localizations  of a ring}. The set $(\Den_l(R), \subseteq )$ is a poset (partially ordered
set). In \cite{larglquot}, it is proved  that the set
$\maxDen_l(R)$ of its maximal elements is a {\em non-empty} set.

$\noindent $

{\it Definition}, \cite{larglquot}. An element $S$ of the set
$\maxDen_l(R)$ is called a {\em maximal left denominator set} of
the ring $R$ and the ring $S^{-1}R$ is called a {\em maximal left
quotient ring} of the ring $R$ or a {\em maximal left localization
ring} of the ring $R$. The intersection
\begin{equation}\label{llradR}
\gll_R:=\llrad (R) := \bigcap_{S\in \maxDen_l(R)} \ass (S)
\end{equation}
is called the {\em left localization radical } of the ring $R$,
\cite{larglquot}.

$\noindent $

 For a ring $R$, there is a canonical exact
sequence 
\begin{equation}\label{llRseq}
0\ra \gll_R \ra R\stackrel{\s }{\ra} \prod_{S\in \maxDen_l(R)}S^{-1}R, \;\; \s := \prod_{S\in \maxDen_l(R)}\, \s_S,
\end{equation}
where $\s_S:R\ra S^{-1}R$, $r\mapsto \frac{r}{1}$.

$\noindent $

{\bf Properties of the maximal left quotient rings of a ring}.
The next theorem describes various properties of the maximal left
quotient rings of a ring, in particular, their groups of units and
their largest left quotient rings.

\begin{theorem}\label{15Nov10}
(\cite{larglquot}.) Let $S\in \maxDen_l(R)$, $A= S^{-1}R$, $A^*$ be
the group of units of the ring $A$; $\ga := \ass (S)$, $\pi_\ga
:R\ra R/ \ga $, $ a\mapsto a+\ga$, and $\s_\ga : R\ra A$, $
r\mapsto \frac{r}{1}$. Then
\begin{enumerate}
\item $S=S_\ga (R)$, $S= \pi_\ga^{-1} (S_0(R/\ga ))$, $ \pi_\ga
(S) = S_0(R/ \ga )$ and $A= S_0( R/\ga )^{-1} R/ \ga = Q_l(R/ \ga
)$. \item  $S_0(A) = A^*$ and $S_0(A) \cap (R/ \ga )= S_0( R/ \ga
)$. \item $S= \s_\ga^{-1}(A^*)$. \item $A^* = \langle \pi_\ga (S)
, \pi_\ga (S)^{-1} \rangle$, i.e. the group of units of the ring
$A$ is generated by the sets $\pi_\ga (S)$ and $\pi_\ga^{-1}(S):=
\{ \pi_\ga (s)^{-1} \, | \, s\in S\}$. \item $A^* = \{ \pi_\ga
(s)^{-1}\pi_\ga ( t) \, |\, s, t\in S\}$. \item $Q_l(A) = A$ and
$\Ass_l(A) = \{ 0\}$.     In particular, if $T\in \Den_l(A, 0)$
then  $T\subseteq A^*$.
\end{enumerate}
\end{theorem}


{\bf A bijection between $\maxDen_l(R)$ and $\maxDen_l(Q_l(R))$}. The next theorem shows that there is a canonical bijection between the maximal left denominator sets of a ring $R$ and its largest left quotient ring $Q_l(R)$. 
\begin{proposition}\label{A8Dec12}
 \cite{Bav-Crit-S-Simp-lQuot} Let $R$ be a ring, $S_l$ be the  largest regular left Ore set of the ring $R$, $Q_l:= S_l^{-1}R$ be the largest left quotient ring of the ring $R$, and $\CC$ be the set of regular elements of the ring $R$. Then
\begin{enumerate}
\item $S_l\subseteq S$ for all $S\in \maxDen_l(R)$. In particular,
$\CC\subseteq S$ for all $S\in  \maxDen_l(R)$ provided $\CC$ is a
left Ore set. \item Either $\maxDen_l(R) = \{ \CC \}$ or,
otherwise, $\CC\not\in\maxDen_l(R)$. \item The map $$
\maxDen_l(R)\ra \maxDen_l(Q_l), \;\; S\mapsto SQ_l^*=\{ c^{-1}s\,
| \, c\in S_l, s\in S\},
$$ is a bijection with the inverse $\CT \mapsto \s^{-1} (\CT )$
where $\s : R\ra Q_l$, $r\mapsto \frac{r}{1}$, and $SQ_l^*$ is the
sub-semigroup of $(Q_l, \cdot )$ generated by the set  $S$ and the
group $Q_l^*$ of units of the ring $Q_l$, and $S^{-1}R= (SQ_l^*)^{-1}Q_l$.
    \item  If $\CC$ is a left Ore set then the map $$ \maxDen_l(R)\ra \maxDen_l(Q), \;\; S\mapsto SQ^*=\{ c^{-1}s\,
| \, c\in \CC, s\in S\}, $$ is a bijection with the inverse $\CT
\mapsto \s^{-1} (\CT )$ where $\s : R\ra Q$, $r\mapsto
\frac{r}{1}$, and $SQ^*$ is the sub-semigroup of $(Q, \cdot )$
generated by the set  $S$ and the group $Q^*$ of units of the ring
$Q$, and $S^{-1}R= (SQ^*)^{-1}Q$.
\end{enumerate}
\end{proposition}

A ring $R$ is called a {\em left Goldie ring} if it satisfies ACC (the {\em ascending chain condition}) for left annihilators  and contains no infinite direct sums of left ideals.

\begin{theorem}\label{A6Jan13}
Let $R$ be a ring such that its left quotient ring $Q$ is a left Noetherian ring. Let $\s :R\ra Q$, $r\mapsto \frac{r}{1}$, $(Q/\gn_Q)^*$ be the group of units of the ring $Q/ \gn_Q$, and $\s':Q/\gn_Q\ra Q(Q/\gn_Q)$, $q+\gn_Q\mapsto \frac{q}{1}+\gn_Q$. Then
\begin{enumerate}
\item $\gn = R\cap \gn_Q$,  $\CC^{-1}\gn = \gn_Q$, $(\CC^{-1}\gn)^i = \CC^{-1}\gn^i$ for all $i\geq 1$, and $\nu = \nu_Q<\infty$ where $\nu $ and $\nu_Q$  are the nilpotency degrees of the prime radicals  $\gn$ and $\gn_Q$,   respectively.
\item
\begin{enumerate}
\item $\CC +\gn \subseteq \CC$.
\item $\tCC \in \Den_l(\bR, 0)$. In particular, $\tCC \subseteq \OCC$.
\item $\tQ := \tCC^{-1}\bR\simeq Q/\gn_Q$ is a semiprime left Noetherian ring.
\item $\gn$ is a nilpotent ideal  of the ring $R$.
\item The $\tQ$-modules $\tCC^{-1} (\gn^i / \gn^{i+1})$, $i=1, \ldots , \nu $,  are   finitely generated  (where $\nu $ is the nilpotency degree of $\gn$).
    \item For each elements $\bc\in \tCC$, the left $\bR$-module $\CN_i/ \CN_i\bc$ is $\tCC$-torsion where $\CN_i:=\gn^i/\gn^{i+1}$.
\end{enumerate}
\item The ring $\bR$ is a semiprime left Goldie ring and its left quotient ring $\bQ := Q(\bR ) \simeq Q(Q/\gn_Q)\simeq Q(\tQ )$ is a semisimple ring.
    \item $1\ra 1+\gn_Q\ra Q^*\stackrel{\pi^*_Q}{\ra}(Q/\gn_Q)^*\ra 1$ is a short exact sequence of group homomorphisms where $\pi_Q:Q\ra Q/ \gn_Q$, $q\mapsto q+\gn_Q$ and $\pi_Q^*:=\pi_Q|_{Q^*}$.
        \item $\CC = \s^{-1}(Q^*)=(\pi_Q\s)^{-1}((Q/\gn_Q)^*)=\pi^{-1} (\ts^{-1}((Q/\gn_Q)^*))$ where $\ts :\bR \ra Q/ \gn_Q$, $\br\mapsto \frac{r}{1}+\gn_Q$, see (\ref{RQbR}).
            \item Let $\dCC:=\CC_{\tQ}$, i.e. $\dCC = \CC_{Q/\gn_Q}$ when we identify the rings $ \tQ$ and $Q/ \gn_Q$ via the isomorphism in statement 2(c). Then $\dCC = \s'^{-1}(Q(Q/\gn_Q)^*)$ and $\OCC = \ts^{-1}(\dCC )=(\s'\ts )^{-1} (Q(Q/\gn_Q)^*)$.
\end{enumerate}
\end{theorem}

{\it Proof}. 1. The prime radical $\gn_Q$ is a nilpotent ideal since the ring $Q$ is a left Noetherian ring. Then  the intersection $R\cap \gn_Q$ is a nilpotent ideal of the ring $R$, hence $R\cap \gn_Q\subseteq \gn$. To establish  the equality $R\cap \gn_Q=\gn$ it suffices to show that the factor ring $R/R\cap \gn_Q$ has no nonzero nilpotent ideals. Suppose that $\bI$ is a nilpotent ideal of the ring $R/ R\cap \gn_Q$ then its preimage $I$ under the epimorphism $R\ra R/ R\cap \gn_Q$, $r\mapsto r+R\cap \gn_Q$, is a nilpotent ideal of the ring $R$ since the ideal $R\cap \gn_Q$ is a nilpotent ideal. We have to show that $\bI =0$. The left ideal $\CC^{-1}I$ of $Q$ is an ideal of the ring $Q$ since the ring $Q$ is a left Noetherian ring. Then $I\CC^{-1}:=\{ ic^{-1}\, | \, i\in I, c\in \CC\} \subseteq \CC^{-1}I$, and so
\begin{equation}\label{CICidN}
(\CC^{-1}I)^i\subseteq \CC^{-1}I^i, \;\;  i\geq 1.
\end{equation}
So, $\CC^{-1}I$ is a nilpotent ideal of $Q$, and so $\CC^{-1}I\subseteq \gn_Q$ and $I\subseteq R\cap \CC^{-1}I\subseteq R\cap \gn_Q$. Therefore, $\bI =0$, as required.
 Clearly,
$$ \CC^{-1} \gn = \CC^{-1} (R\cap \gn_Q)= \CC^{-1} R\cap \CC^{-1}\gn_Q= Q\cap \gn_Q=\gn_Q,$$and $(\CC^{-1}\gn )^i = \CC^{-1}\gn^i$ for $i\geq 1$. In particular, $\nu = \nu_Q<\infty$.

2(a). Let $c\in \CC$ and $n\in \gn$. Then the element $c^{-1}n\in \CC^{-1} \gn = \gn_Q$ (statement 1)  is a nilpotent element of the ring $Q$ and so the element  $1+c^{-1}n$ is a unit of the ring $Q$. Now,
$$c+n =c(1+c^{-1}n)\in \CC.$$
2(b,c) Since $ \gn = R\cap \gn_Q$ (statement 1), there is a commutative diagram  of ring homomorphisms
\begin{equation}\label{RQbR}
\xymatrix{Q\ar[r]^{\pi_Q} & Q/\gn_Q\\
R\ar[r]^\pi\ar[u]^\s &\bR\ar[u]_{\ts}}\\
\end{equation}
where the horizontal maps are natural epimorphisms  and the vertical maps are natural monomorphisms (where $\ts ( \br ) := \frac{\br}{1}=\frac{r}{1}+\gn_Q$). Since
$$Q/\gn_Q=\{ \pi (c)^{-1} \br \, | \, c\in \CC , \br\in \bR \},$$ we see that $\tCC = \pi (\CC )\in \Den_l(\bR , 0)$ and $Q/ \gn_Q\simeq \tCC^{-1}\bR$ is a semiprime left Noetherian ring.

2(d)  Statement 1.

2(e) By the statement (c), $\tQ \simeq Q/ \gn_Q$ is a left Noetherian ring. Hence, the $Q/ \gn_Q$-modules
$$\gn_Q^i / \gn_Q^{i+1} \simeq (\CC^{-1} \gn )^i / (\CC^{-1} \gn )^{i+1} \stackrel{{\rm st.1}}{\simeq} \CC^{-1}\gn^i / \CC^{-1}\gn^{i+1}\simeq \tCC^{-1} (\gn^i/ \gn^{i+1})$$ are finitely generated where $i=1, \ldots , \nu$.

2(f) For each $i=1, \ldots , \nu$, the left $Q$-module/$\tQ$-module
 $\CN_i/ \CN_i\bc$ is  $\tCC$-torsion since
\begin{eqnarray*}
\tCC^{-1}(\CN_i/ \CN_i\bc)&=& \tCC^{-1} (\gn^i/(\gn^ic+\gn^{i+1}))=\CC^{-1}(\gn^i/(\gn^ic+\gn^{i+1}))\\
&=& \CC^{-1}\gn^i /(\CC^{-1}\gn^ic+\CC^{-1} \gn^{i+1})=\CC^{-1}\gn^i/\CC^{-1}\gn^i=0.
\end{eqnarray*}

4. Statement 4 follows from the fact that $\gn_Q$ is a nilpotent ideal of the ring $Q$.

3. By statement 2(c), the ring $\tQ\simeq Q/ \gn_Q$ is a semiprime left Noetherian ring. In particular, it is a semiprime left Goldie ring and, by Goldie's Theorem, its left quotient ring $Q(\tQ )\simeq Q(Q/ \gn_Q)$ is a semisimple ring. Since $\bR \subseteq Q/ \gn_Q \simeq \tCC^{-1}\bR$ (statement 2(c)) we have $Q(\bR )=Q(\tCC^{-1}\bR )$ is a semisimple ring. By Goldie's Theorem, the ring $\bR$ is a semiprime left Goldie ring. So, we can extend the commutative diagram (\ref{RQbR}) to the commutative diagram (which is used in the proof of statement 6)

\begin{equation}\label{RQbR1}
\xymatrix{Q\ar[r]^{\pi_Q} & Q/\gn_Q\ar[r]^{\s'}        &   Q(Q/\gn_Q)\ar[d]^{\simeq}\\
R\ar[r]^\pi\ar[u]^\s      & \bR\ar[u]_{\ts}\ar[r]^{\overline{\s}} &    \bQ }
\end{equation}
where the maps $\s'$ and $\overline{\s}$ are monomorphisms, $\s'(q+\gn_Q) = \frac{q+\gn_Q}{1}$ and $\overline{\s}(\br ) = \frac{\br}{1}$.

5. By Theorem  \ref{4Jul10}.(1),  $\CC =\s^{-1}(Q^*)$. By statement 4, $Q^* = \pi_Q^{-1} ((Q/ \gn_Q)^*)$. Then, in view of the commutative diagram  (\ref{RQbR}),
\begin{eqnarray*}
 \CC &=& \s^{-1}\pi_Q^{-1} ((Q/ \gn_Q)^*)=(\pi_Q\s)^{-1}((Q/ \gn_Q)^*)=(\ts \pi)^{-1}((Q/ \gn_Q)^*)\\
 &=&\pi^{-1}(\ts^{-1}((Q/ \gn_Q)^*))).
\end{eqnarray*}

6. By Theorem \ref{4Jul10}.(1),  $\dCC = \s'^{-1} (Q(Q/ \gn_Q)^*)$ and $\OCC = \overline{\s}^{-1}(\bQ^*)=\overline{\s}^{-1}(Q(Q/ \gn_Q)^*)$ (statement 3).  Thus, the commutativity of the second square in the diagram (\ref{RQbR1}) yields,
$$ \OCC = (\s'\ts )^{-1}(Q(Q/ \gn_Q)^*)=\ts^{-1} (\dCC ).\;\;\;\; \Box $$

$\noindent $

The next corollary is a criterion for a ring to have a left Noetherian left quotient  ring $Q$ such that the factor ring $Q/\gn_Q$ is a semisimple ring (or $\tQ$ is a semisimple ring; or $\tCC = \OCC$; or $\CC = \pi^{-1}(\OCC )$).

\begin{corollary}\label{a8Jan13}
Let $R$ be a ring such that its left quotient ring $Q$ is a left Noetherian ring, we keep the notation of Theorem \ref{A6Jan13} and its proof. The following statements are equivalent (recall that $\tQ=Q/\gn_Q$, Theorem \ref{A6Jan13}.(2c)).
\begin{enumerate}
\item $\tQ$ is a semisimple ring.
\item $\tQ = Q(\tQ )$.
\item $\OCC = \ts^{-1}(\tQ^*)$.
\item $\CC = \pi^{-1}(\OCC )$.
\item $\tCC = \OCC$.
\item $\tQ$ is a left Artinian ring.
\end{enumerate}
\end{corollary}

{\it Proof}. $(1\Rightarrow 2)$ Trivial.

$(2\Rightarrow 3)$ If $\tQ = Q(\tQ )$, i.e. the map $\s'$ in (\ref{RQbR1}) is an isomorphism, then the rings $ \tQ \simeq Q/ \gn_Q$ and $\bQ$ are isomorphic, see the commutative diagram  (\ref{RQbR1}). Now, $\OCC =\overline{\s}^{-1} (\bQ^*) = \ts^{-1} (\tQ^*)$ where the first equality holds by Theorem \ref{4Jul10}.(1).

$(3\Rightarrow 4)$ By Theorem \ref{A6Jan13}.(5). $\CC =\pi^{-1} (\ts^{-1}(\tQ^*))=\pi^{-1} (\OCC )$.

$(4\Rightarrow 5)$ $\tCC =\pi (\CC ) = \pi (\pi^{-1} (\OCC )) = \OCC$ since the map $\pi$ is an epimorphism.

$(5\Rightarrow 1)$ If $\tCC = \OCC$ then by (\ref{RQbR1}),  $\tQ \simeq  \bQ$ is a semisimple ring, by Theorem \ref{A6Jan13}.(3).

$(1\Rightarrow 6)$ Trivial.

$(6\Rightarrow 1)$ This implication follows from Theorem \ref{A6Jan13}.(2c). $\Box $



\section{Proof of Theorem \ref{6Jan13} and a Second Criterion (via the associated graded ring)}\label{PTHA6Jan13}
 The aim of this section is to give  proofs of Theorem \ref{6Jan13} and  Theorem \ref{27Jan13} which are criteria for a ring $R$ to have a left Noetherian left quotient ring.

{\bf Proof of Theorem \ref{4Oct14}}. The implications $(1\Leftarrow 2)$ and  $(4\Leftarrow 5)$ are  obvious and the equivalence $(2\Leftrightarrow 3)$ is Goldie's Theorem.

$(1\Rightarrow 2)$ Since $R$ is a semiprime ring then so is the ring $Q=Q(R)$. In particular, the ring $Q$ is a semiprime left  Goldie ring, hence, by Goldie's Theorem, $Q(Q(R))= Q(R)$ is a semisimple ring.

$(2\Leftrightarrow 5)$ {\cite[Theorem 2.9]{larglquot}}.

$(4\Rightarrow 2)$ Since $R$ is a semiprime ring then so is the ring $Q_l(R)$ (since the ring $Q_l(R)$ is left Noetherian). The ring $Q_l(R)$ is a left Noetherian ring, hence $Q_l(R)$ is a semiprime left Goldie ring. By Goldie's Theorem, $Q(Q_l(R))$ is a semisimple ring.
  Then, by {\cite[Theorem 2.9]{larglquot}}, $Q(Q_l(R))= Q_l(Q_l(R))$ is a semisimple ring. By {\cite[Theorem 2.8.(4)]{larglquot}},  $Q_l(Q_l(R))= Q_l(R)$, and so $Q_l(R)$ is a semisimple ring. Then, by  {\cite[Theorem 2.9]{larglquot}}, $Q(R) = Q_l(R)$, i.e. $Q(R)$ is a semisimple ring. $\Box$

{\bf Proof of Theorem \ref{6Jan13}}. $(1\Rightarrow 2)$ Theorem \ref{A6Jan13}.(1,2).

$(1\Leftarrow 2)$ (i) $\CC \in \Ore_l(R)$: We have to show that for given  elements $c\in \CC$ and $r\in R$ there are elements $c'\in \CC$ and $r'\in R$ such that $c'r=r'c$.  We can assume that $r\neq 0$ since otherwise take $c'=1$ and $r'=0$. To prove this fact  we use a downward induction on the {\em degree} of the element $r\neq 0$:
$$ \deg (r) :=\max \{ i\, | \, r\in \gn^i \;\; {\rm where}\;\;  0\leq i\leq \nu\}.$$
Suppose that $\deg (r) = \nu$. By the condition (f), the $\bR$-module $\gn^\nu/ \gn^\nu \bc $ is $\tCC$-torsion, hence $\bc'r=r'\bc$ for some elements $\bc'=c+\gn \in \tCC$, $c\in \CC$ and $r'\in \gn^\nu$. Then $c'r=r'c$.

Suppose that $i=\deg (r) <\nu$, and the result is true for all elements $c\in \CC$ and $r\in R$ with $\deg (r)>i$.

Suppose that $i=0$, i.e. $r\in R\backslash \gn$. Since $\tCC \in \Den_l(\bR , 0)$ (by the conditions (a) and (b)), $\bc_1\br = \br_1\bc$ for some elements $c_1\in \CC$ and $r_1\in R$. The difference $a:= c_1r-r_1c$ belongs to the ideal $\gn$, and so $\deg (a)>0$.  By induction, $c_2a=bc$ for some elements $c_2\in \CC$ and $b\in R$. Then
$$ c_2c_1r=(c_2r_1+b)c,$$ and it suffices to take $c'=c_2c_1$ and $r' = c_2r_1+b$.

Suppose that $i>0$. By the condition (f),  the left $\bR$-module $\CN_i/ \CN_i\bc=\gn^i/(\gn^ic+\gn^{i+1})$ is $\tCC$-torsion. Therefore, $sr=xc+y$ for some elements $s\in \CC$, $x\in \gn^i$ and $y\in \gn^{i+1}$. Since $\deg (y) \geq i+1$, by induction, there are elements $t\in \CC$ and $z\in R$ such that $ty=zc$. Therefore, $$tsr=(tx+z)c$$ and it suffices to take $c'=ts$ and $r'= tx+z$.

(ii) {\em $\CC^{-1}\gn$ is an ideal of the ring $Q$ such that $Q/ \CC^{-1}\gn \simeq \tQ$}:  By (i), the left quotient ring $Q= \CC^{-1}R$ exists. Let $\s : R\ra Q$, $r\mapsto \frac{r}{1}$. By the universal property of left localization, there is a ring homomorphism $\pi_Q:Q\ra \tQ$, $ c^{-1}r\mapsto \bc^{-1}\br$, where $\bc = c+\gn$ and $\br = r+\gn$, and we have the commutative diagram of ring homomorphisms

\begin{equation}\label{RQbR3}
\xymatrix{Q\ar[r]^{\pi_Q} & \tQ\\
R\ar[r]^\pi\ar[u]^\s &\bR\ar[u]_{\ts}}\\
\end{equation}

where $\ts : \br \mapsto \frac{\br}{1}$ is a monomorphism and $\pi_Q$ is an epimorphism (by the very definition of $\pi_Q$). Applying the exact functor $\CC^{-1} (-)$ to the short exact sequence of $R$-modules $0\ra \gn \ra R\stackrel{\pi}{\ra} \bR \ra 0$ we obtain
the short exact sequence of $Q$-modules
$$ 0\ra \CC^{-1}\gn \ra Q\stackrel{\pi_Q}{\ra} \CC^{-1}\bR =\tCC^{-1}\bR =\tQ \ra 0.$$
 Therefore, $\ker (\pi_Q) = \CC^{-1} \gn$ is an ideal of $Q$ (since $\pi_Q$ is a ring homomorphism)  such that $Q/\CC^{-1}\gn\simeq \tQ$.

(iii) $(\CC^{-1}\gn )^i = \CC^{-1}\gn^i$ {\em for} $i\geq 1$: This follows from (ii).

(iv) {\em The ring $Q$ is a left Noetherian ring}: By localizing the descending chain of ideals of the ring $R$:
$$ R\supset \gn \supset \gn^2\supset \cdots \supset \gn^i \supset \cdots \supset \gn^\nu \supset \gn^{\nu +1} =0$$
we obtain the descending chain of ideals (by (iii)) of the ring $Q$:
\begin{equation}\label{QCCf}
Q\supset \CC^{-1}\gn \supset \CC^{-1}\gn^2\supset \cdots \supset \CC^{-1}\gn^i \supset \cdots \supset \CC^{-1}\gn^\nu \supset \CC^{-1}\gn^{\nu +1} =0.
\end{equation}
By (ii), $\tQ\simeq Q/ \CC^{-1}\gn$ is a left Noetherian $Q$-module since the ring $\tQ$ is a left Noetherian ring, by the condition (c). For each $i=1, \ldots , \nu$, the left $Q$-module, $\CC^{-1}\gn^i / \CC^{-1} \gn^{i+1}\simeq (\CC^{-1}\gn )^i / (\CC^{-1} \gn )^{i+1}$ is a $Q/ \CC^{-1} \gn = \tQ$-module, by (ii). The $\tQ$-modules
$$\CC^{-1}\gn^i / \CC^{-1} \gn^{i+1}\simeq \CC^{-1}(\gn^i / \gn^{i+1})\simeq \tCC^{-1}(\gn^i / \gn^{i+1})=\tCC^{-1}\CN_i$$ are finitely generated, by the condition (e), hence Noetherian since the ring $\tQ$ is a left Noetherian. Since all the factors of the finite filtration (\ref{QCCf}) are Noetherian $\tQ$-modules/$Q$-modules, the ring $Q$ is a left Noetherian ring. $\Box $

$\noindent $

{\bf The minimal primes of the rings $R$, $Q(R)$, $\bQ$ and $\tQ$}. For a ring $R$ such that $Q(R)$ is a left Noetherian ring, the next corollary shows that the localizations of $R$ at the maximal left denominator sets are left Noetherian rings and there are natural bijections between the sets of minimal primes of the rings  $R$, $Q(R)$, $\bQ$ and $\tQ$.

\begin{corollary}\label{a7Jul15}
Let $R$ be a ring such that $Q(R)$ is a left Noetherian ring. Then 
\begin{enumerate}
\item For every $S\in \maxDen_l(R)$, the ring $S^{-1}R$ is a left Noetherian ring.
\item
\begin{enumerate}
\item The map $\Min(R)\ra \Min (\bQ )$, $\gp\mapsto \OCC^{-1}(\gp /\gn )$, is  a bijection with the inverse $\gq\mapsto (\overline{\s}\pi )^{-1}(\gq )$ where the maps $\overline{\s}$ and $\pi$ are defined in (\ref{piRRQ}). 
\item The map $\Min (R)\ra \min (\tQ )$, $\gp\mapsto \tCC^{-1}(\gp / \gn )$, is a bijection with the inverse $\gq\mapsto \tau^{-1}(\gq )$ where $\tau : R\ra \tQ$, $r\mapsto\frac{r+\gp}{1}$.
\item The map $\Min (R)\ra \Min (Q)$, $\gp\mapsto \CC^{-1}\gp$, is a bijection with the inverse $\gq \mapsto \gq \cap R$.
\end{enumerate}
\end{enumerate}
\end{corollary}

{\it Proof}. 1. By Proposition \ref{A8Dec12}.(1), $\CC =\CC_R\subseteq S$ for all $S\in \maxDen_l(R)$. The ring $Q:=\CC^{-1}R$ is a left Noetherian ring and the ring $S^{-1}R\simeq (SQ^*)^{-1}Q$ is a left localization of the ring $Q$ (by Proposition \ref{A8Dec12}.(3), the map $\maxDen_l(R)\ra \maxDen_l(Q)$, $S\mapsto SQ^*$, is a bijection with $S^{-1}R\simeq (SQ^*)^{-1}Q$). Therefore, the ring $S^{-1}R$ is a left Noetherian ring since the ring $Q$ is so. 

2(a). The map $\Min (R)\ra \Min (\bR )$, $\gp\mapsto \gp /\gn$, is a bijection with the inverse $\gp'\mapsto \pi^{-1}(\gp')$ where $\pi :R\ra \bR$, $r\mapsto r+\gp$. The ring $\bR$ is a semiprime left  Goldie ring such that $\bQ\simeq Q(Q/\gn_Q)$ is a semisimple ring. Hence, the map $\min (\bR)\ra \Min (\bQ)$, $ \gp'\mapsto \OCC^{-1}\gp'$, is a bijection with the inverse $\gq \mapsto \overline{\s}^{-1}(\gq )$ where 
 $\overline{\s}: \bR\ra \bQ$, $\br\mapsto \frac{\br}{1}$. Now, the statement (a) follows.
 
 (b) The ring $\tQ$ is a semiprime left Goldie ring (by Theorem \ref{A6Jan13}.(2c)) and $Q(\tQ )\simeq \tQ$ is a semisimple ring (Theorem \ref{A6Jan13}.(3)). So, the map $\Min (\tQ )\ra \Min (\bQ )$, $ P\mapsto \bQ\t_{\tQ}P$, is a bijection with the inverse $P'\mapsto P'\cap \tQ$. Now, the statement (b) follows from the statement (a).
 
 (c) By Theorem \ref{A6Jan13}.(2c), $\tQ\simeq Q/\gn_Q$ and the map $\Min (Q)\ra \Min (Q/\gn_Q)$, $P\mapsto P/\gn_Q$, is a bijection. Now, the statement (c) follows from the statement (b). $\Box $


\begin{lemma}\label{a27Jan13}
Let $G$ be a monoid and $e\in G$ be its neutral element, $A=\bigoplus_{g\in G} A_g$ be a $G$-graded ring, $1\in A_e$, $S\in \Den_l (A)$ and $\ga = \ass (S)$. If $S\subseteq A_e$ then the ring $S^{-1}A=\bigoplus_{g\in G}(S^{-1}A)_g$ is a $G$-graded ring where $(S^{-1}A)_g=S^{-1}A_g:=\{ s^{-1}a_g\, | \, s\in S, a_g\in A_g\}$, $S\in \Den_l(A_e)$ and $\ga =\bigoplus_{g\in G} \ga_g$ is a $G$-graded ideal of the ring $A$, i.e. $\ga_g=\ga \cap A_g$ for all $g\in G$.
\end{lemma}

{\it Proof}. Straightforward.  $\Box $


Suppose that a ring $R$ has a left Noetherian left quotient ring $Q$. By Theorem \ref{A6Jan13}, the associated graded ring $\grQ := Q/\gn_Q\oplus \gn_Q/\gn_Q^2\oplus \cdots $ is equal to
$$\grQ =\tQ\oplus \tCC^{-1}(\gn / \gn^2)\oplus \cdots \oplus \tCC^{-1}(\gn^\nu / \gn^{\nu +1}).$$

{\bf  Proof of Theorem \ref{27Jan13}}. $(1\Rightarrow 2)$ Suppose that the ring $Q$ is a left Noetherian ring, and so the conditions of Theorem \ref{6Jan13} and Theorem \ref{A6Jan13} hold. In particular, $\tCC\subseteq \OCC$ and $\gn$ is a nilpotent ideal.

(i) $\tCC\in \Ore_l(\grR )$: It suffices to show that for given elements $\bc = c+\gn\in \tCC$ and $ r+\gn^{i+1}\in \gn^i/\gn^{i+1}$ where $c\in \CC$, $r\in \gn^i$ and $i=0,1, \ldots , \nu$, there are elements  $\bc' = c'+\gn\in \tCC$ and $ r'+\gn^{i+1}\in \gn^i/\gn^{i+1}$ where $c'\in \CC$, $r'\in \gn^i$ such that $\bc'(r+\gn^{i+1}) = (r'+\gn^{i+1}) \bc$.

The case $i=0$ is obvious, by Theorem \ref{6Jan13}.(2b). So, we can assume that $i\geq 1$. By Theorem \ref{A6Jan13}.(1), $$\gn^i\CC^{-1}:=\{ nc^{-1}\, | \, n\in \gn^i , c\in \CC\} \subseteq \CC^{-1}\gn^i.$$ So, $rc^{-1} = c'^{-1}r'$ for some elements $c'\in \CC$ and $r'\in \gn^i$, hence $ c'r=r'c$. This equality implies the required one.

(ii) $\tCC\in \Den_l(\grR )$: We have to show that if $\br \bc =0$ for some elements $\br = r+\gn^{i+1}\in \gn^i/\gn^{i+1}$ and $\bc =c+\gn$ where $r\in \gn^i$, $c\in \CC$ and $i=0,1,\ldots , \nu$, then $\bc_1\br =0$ for some element $\bc_1=c_1+\gn$ where $c_1\in \CC$. The case $i=0$ follows from Theorem \ref{A6Jan13}.(2b). Suppose that $i\geq 1$. Then $n:=rc\in \gn^{i+1}$, and so $r=nc^{-1} = c_1^{-1}n_1$ for some elements $c_1\in \CC$ and $ n_1\in \gn^{i+1}$, by Theorem \ref{A6Jan13}.(2c). Hence, $\bc_1\br =0$, as required. 

(iii) $\grQ\simeq \tCC^{-1}\grR$: By Theorem \ref{A6Jan13}.(1), $\gn_Q=\CC^{-1}\gn$. By Theorem \ref{A6Jan13}.(2c), $Q/\gn_Q\simeq \tQ$. Now, using Theorem \ref{A6Jan13}.(1,2), we have
\begin{eqnarray*}
\grQ &= &  \tQ\oplus \cdots \oplus \gn_Q^i/\gn_Q^{i+1}\oplus\cdots = \tQ\oplus \cdots \oplus \CC^{-1}\gn^i/\CC^{-1}\gn^{i+1}\oplus\cdots \\
 & =&  \tQ\oplus \cdots \oplus \CC^{-1}(\gn^i/\gn^{i+1})\oplus\cdots =\tQ\oplus \cdots \oplus \tCC^{-1}(\gn^i/\gn^{i+1})\oplus\cdots\\
 &\simeq &\tCC^{-1}\grR.
\end{eqnarray*}
(iv) $\tCC^{-1}\grR$ {\em is a left Noetherian ring}: The ring $\tQ$ is a left Noetherian ring (Theorem \ref{A6Jan13}.(2c)) and the left $\tQ$-modules $\tCC^{-1}(\gn^i / \gn^{i+1})\simeq (\CC^{-1}\gn )^i / (\CC^{-1}\gn )^{i+1}$ are finitely generated where $i=1, \ldots , \nu$ (since $Q$ is a left Noetherian ring). Therefore, the left $\tQ$-module $\grQ= \tQ\oplus \cdots \oplus \tCC^{-1} (\gn^i/\gn^{i+1})\oplus \cdots \oplus \tCC^{-1} (\gn^\nu/\gn^{\nu+1})$ is finitely generated, hence Noetherian. Since $\tQ \subseteq \grQ$, the ring $\grQ$ is a left Noetherian ring.

$(1\Leftarrow 2)$ It suffices to show that the conditions (a)-(f) of Theorem \ref{6Jan13}.(2) hold. The conditions (a) and (d) are given. The set $\tCC$ is a left denominator set of the $\N$-graded ring $\grR$ such that $\tCC\subseteq \bR$. By Lemma \ref{a27Jan13}, the ring $\tCC^{-1} \grR= \tQ\oplus \cdots \oplus \tCC^{-1}\CN_i\oplus\cdots$ is an $\N$-graded ring and $\tCC\in \Den_l(\bR , 0)$, and so the condition (b) holds.

The ring $\tQ$ is a factor ring of the left Noetherian ring $\tCC^{-1}\grR$, hence $\tQ$ is a left Noetherian ring, i.e. the condition (c) holds.

The ring $\tCC^{-1}\grR$ is left Noetherian, hence the $\tQ$-modules $\tCC^{-1}\CN_i$ are finitely generated, i.e. the condition (e) holds.

The ring  $\tCC^{-1}\grR$ is an $\N$-graded ring. In particular, $\tCC^{-1} \CN_i\tQ\subseteq \tCC^{-1}\CN_i$ for all $i$. Therefore, $\CN_i\tCC^{-1}=\{ nc^{-1}\, | \, n\in \CN_i, c\in \tCC\}\subseteq \tCC^{-1}\CN_i$, i.e. the condition (f) holds. $\Box $


\begin{corollary}\label{x28Jan13}
Let $R$ be a  ring with a left Noetherian left quotient ring $Q$ and $\tga :=\ass_{\grR }(\tCC )$. Then the largest left quotient ring $Q_l(\grR / \tga )$ of the ring $\grR / \tga $ is a left Noetherian ring.
\end{corollary}

{\it Proof}. By Theorem \ref{27Jan13}.(2), the ring $\tCC^{-1} \grR $ is a left Noetherian ring. The ring $Q_l(\grR / \tga )$  is a left localization  of the ring  $\tCC^{-1} \grR $, hence is a left  Noetherian ring. $\Box $


\begin{proposition}\label{a28Jan13}
Let $R$ be a ring and $S\in \Den_l(R)$. Then the ring $R$ is a left Noetherian ring iff the ring $S^{-1}R$ is a left Noetherian ring  and for each left ideal $I$ of $R$ the left $R$-module $\tor_S(R/I)$ is  left Noetherian.
\end{proposition}

{\it Proof}. $(\Rightarrow )$ Trivial.

$(\Leftarrow )$ Let $I$ be a left ideal of the ring $R$. We have to show that $I$ is a finitely generated left $R$-module. The ring $S^{-1}R$ is a left Noetherian ring. Then its left ideal $S^{-1}I$ is a finitely generated $S^{-1}R$-module. We can find elements $u_1, \ldots , u_n\in I$ such that $S^{-1}I=\sum_{i=1}^nS^{-1}Ru_i/1$. Let $I':=\sum_{i=1}^nRu_i$. Then $I'\subseteq I$ and $I/I'$ is a submodule of the left Noetherian $R$-module $\tor_S(R/I')$ and as a result the left $R$-module $I/I'$ is finitely generated. This fact implies that the left ideal $I$ is finitely generated.  $\Box $

$\noindent $

Let $R$ be a ring, $S\in \Ore_l(R)$ and $M$ be a left $R$-module. Let $\ker (S, M) = \{ \ker (s_M)\, | \,  s\in S\}$ where $s_M: M\ra M$, $ m\mapsto sm$. The set $(\ker (S, M) , \subseteq )$ is a poset. Let $ \max\ker (S, M)$ be the set of its maximal elements and let $\max (S, M) := \{ s\in S\, | \, \ker (s_M) \in  \max\ker (S, M)\}$.

\begin{proposition}\label{a5Oct14}
Let $R$ be a ring, $S\in \Den_l(R)$, $M$ be a left $R$-module such that $ \max\ker (S, M) \neq \emptyset$. Then $\max\ker (S, M)=\{ \tor_S(M)\}$, i.e. $\ker (s_M)= \tor_S(M)$ for all elements $s\in \max (S, M)$. In particular, $\ga := \ann_R(\tor_S(M))\neq 0$ and $\emptyset \neq \max (S, M) \subseteq S\cap \ga$.
\end{proposition}

{\it Proof}. Suppose that $\ker (s_M) \neq \tor_S(M)$ for some $s\in \max (S, M)$, we seek a contradiction. Fix $m\in  \tor_S(M)\backslash \ker (s_M)$ and $t\in S$ with $tm =0$. Since $S\in \Ore_l(R)$, we have $St\cap Rs\neq \emptyset$, i.e. $s' := s_1t\in St\cap Rs$ for some $s_1\in S$. Clearly, $s'\in S$ and $\ker (s_M') \supseteq \ker (t_M) + \ker(s_M) \varsupsetneqq \ker (s_M)$, a contradiction. $\Box $


\begin{proposition}\label{b5Oct14}
Let $R$ be a ring, $\gp_1, \ldots , \gp_n$ be prime ideals of the ring $R$ such that the rings $R/ \gp_i$ $(i=1, \ldots , n)$ are left Goldie rings. Let $M$ be a left $R$-module such that $M= M_n\supset M_{n-1}\supset \cdots \supset M_1\supset M_0=0$ is a chain of its submodules such that $\gp_i$ is  maximal
 among left annihilators of non-zero submodules of $M/ M_{i-1}$, $M_i= \ann_{M/M_{i-1}}(\gp_i) := \{ m\in M/ M_{i-1} \, | \, \gp_im=0\}$, and $\max\ker (\CC_{R/ \gp_i}, R/ M_{i-1})\neq \emptyset$ for $i=1, \ldots , n$. Then $\bigcap_{i=1}^n \CC (\gp_i) \subseteq {}'\CC_M(0):=\{ r\in R\, | \, r_M$ is an injection$\}$.
 \end{proposition}

{\it Proof}. Let $c\in \bigcap_{i=1}^n \CC (\gp_i)$. We have to show that $\ker (c_M)=0$. Suppose that this is not true, i.e. $cm=0$ for some $0\neq m\in M$, we seek a contraction. Then $m\in M_i\backslash M_{i-1}$ for some $i$. Then $0\neq \overline{m}:= m+M_{i-1} \in \tau_i:= \tor_{\CC_{R/\gp_i}}(M_i/ M_{i-1})$ 
 since $\overline{c}:= c+\gp_i\in \CC_{R/\gp_i}$. By Proposition \ref{a5Oct14}, $\overline{s}\tau_i=0$ for some $\overline{s}\in  \CC_{R/\gp_i}$ (since $\max\ker ( \CC_{R/\gp_i},R/M_{i-1})\neq \emptyset$), i.e.  $s\tau_i=0$ for some $s\in \CC (\gp_i)$. This means that $\ann_R(\tau_i) \varsupsetneqq \gp_i$, a contradiction. $\Box$


\section{Finiteness of $\maxDen_l(R)$ when $Q(R)$ is a left Noetherian ring}\label{FMDQRLN}

 The aim of this section is to give a proof of Theorem \ref{aA5Oct14}.

Let $R$ be a ring. Let $S$ and $T$ be submonoids of the multiplicative monoid $(R, \cdot )$. We denote by $ST$ the {\em submonoid}  of $(R, \cdot )$ generated by $S$ and $T$. This notation should not be confused with the product of two sets  which is {\em not} used in this paper. The next result is a criterion for the set $ST$ to be a left Ore (denominator) set, it is used at the final stage in the proof of Theorem \ref{aA5Oct14}.

\begin{lemma}\label{a14Sep13}
(\cite{Bav-stronglquot}.)
\begin{enumerate}
\item Let $S,T\in \Ore_l(R)$. If $0\not\in ST$ then $ST \in \Ore_l (R)$.  \item Let $S,T\in \Den_l(R)$. If $0\not\in ST$ then $ST \in \Den_l (R)$.
\end{enumerate}
\end{lemma}


{\bf  Proof of Theorem \ref{aA5Oct14}}. Let $S\in \maxDen_l(R)$ and $\pi : R\ra \bR$, $r\mapsto \br :=r+\gn$.

(i) $\bS:= \pi (S)\in \Ore_l(\bR )$: This obvious since $S\cap \gn = \emptyset$.

(ii) $\CC \subseteq S$ (Proposition \ref{A8Dec12}.(1)).

(iii) $S^{-1}R$ {\em is a left Noetherian ring}:  By (ii), $S^{-1}R$ is a left localization of the left Noetherian ring $Q$, hence $S^{-1}R$ is a left Noetherian ring.

(iv) $S^{-1}\gn$ {\em is an ideal of} $S^{-1}R$ {\em such that} $(S^{-1}\gn )^i = S^{-1} \gn^i$ {\em for all} $i\geq 1$: By (iii), $S^{-1}\gn$ is an ideal of the ring $S^{-1}R$. In particular, $\gn S^{-1} \subseteq S^{-1} \gn$, hence $(S^{-1}\gn )^i = S^{-1} \gn^i$  for all $i\geq 1$.

(v) $\bS \in \Den_l(\bR )$: In view of (i), we have to show that if  $\br\bs =0$ for some $\br \in \bR$ and $\bs\in \bS$ then $\bt \br =0$ for some element $\bt \in \bS$. The equality
 $\br\bs =0$ means that $n:=rs\in \gn$. Then $S^{-1} R \ni \frac{r}{1}=ns^{-1} = s_1^{-1}n_1$ for some $s_1\in S$ and $ n_1\in \gn$, (by (iv)). Hence, $s_2s_1r= s_2n_1\in \gn$ for some $s_2\in S$, and so $\bt \br =0$ where $\bt = \bs_2\bs_1\in \bS$.

 The ring $\bR$ is a semiprime left Goldie ring (Theorem \ref{A6Jan13}.(3)), $\bQ= \prod_{i=1}^s\bQ_i$ and $\bQ_i$ are simple Artinian rings. By  \cite[Theorem 4.1]{Bav-Crit-S-Simp-lQuot}, $|\maxDen_l(\bR )|=s$. Let $\maxDen_l(\bR ) = \{ T_1, \ldots , T_s\}$.

(vi) $|\maxDen_l(R )|\leq s$: By (v), $\bS\subseteq T_i$ for some $i$. Then $S$ is the largest element with respect to inclusion of the set $\{ S'\in \Den_l(R)\, | \, \pi (S')\subseteq T_i\}$ (as $0\in SS_1$ for all {\em distinct} $S,S_1\in \maxDen_l(R)$, by Lemma \ref{a14Sep13}.(2), where $SS_1$ is the multiplicative submonoid of $R$ generated by $S$ and $S_1$). Hence, $|\maxDen_l(R )|\leq s$. $\Box $

$\noindent $

Let $R$ be a ring such that its left quotient ring $Q(R)$ is a left Noetherian ring. Define the ring homomorphism
\begin{equation}\label{piRRQ}
  p_i: R\stackrel{\pi}{\ra} \bR \stackrel{\overline{\s} }{\ra}\bQ=\prod_{j=1}^s\bQ_i\stackrel{\bp_i}{\ra} \bQ_i. 
\end{equation}
Let $\overline{\s}_i:=\bp_i\overline{\s}:\bR \ra \bQ_i$. By \cite[Theorem 4.1]{Bav-Crit-S-Simp-lQuot}, $\maxDen_l(\bR )=\{ T_1, \ldots , T_s\}$ where $T_i:= \overline{\s}_i^{-1}(\bQ_i^*)=\{ \br\in \bR\, | \, \overline{\s}_i(\br ) \in \bQ_i^*\}$. 

$\noindent$

{\bf Proof of Corollary \ref{xaA5Oct14}}. 1. The set $D_i$ is a non-empty set since $R^*\subseteq D_i$. By Lemma \ref{a14Sep13}.(2), the set $S_i':=\bigcup_{S'\in D_i}S'$ is the largest element of the set $D_i$. In the proof of Theorem \ref{aA5Oct14}, we proved that $\maxDen_l(R)\subseteq \CM$ where $\CM$ is the set of maximal elements of the set $\{ S_1',\ldots , S_s'\}$ (see the statement (vi)). The reverse inclusion is obvious.

2. By Theorem \ref{A6Jan13}.(2b), $\tCC \in \Den_l(\bR , 0)$ and $\tCC \subseteq \OCC$. By Theorem \ref{A6Jan13}.(3), the ring $\bR$ is a
 semiprime left Goldie ring and $\bQ$ is a semisimple ring. By Theorem \ref{4Jul10}.(1), $\OCC=\bR\cap \bQ^*$. Hence, $\CC\subseteq S_i'$ for all $i=1,\ldots , s$.

3. Statement 3 follows from statement 2: The ring ${S_i'}^{-1}R$ is a localization of the left Noetherian ring $\CC^{-1}R$ (see statement 2). Hence it is a left Noetherian ring. $\Box$

$\noindent$

{\bf Criterion for $|\maxDen_l(R) |= |\maxDen_l(\bR ) |$}. Let $R$ be a ring such that its left quotient ring $Q(R)$ is a left Noetherian ring.   In general, $|\maxDen_l(R) |\leq  |\maxDen_l(\bR ) |$, Theorem \ref{aA5Oct14}. The next theorem is a criterion for  $|\maxDen_l(R) |= |\maxDen_l(\bR ) |$.

\begin{theorem}\label{A23Nov14}
 Let $R$ be a ring such that its left quotient ring $Q(R)$ is a left Noetherian ring. Then $|\maxDen_l(R) |= |\maxDen_l(\bR ) |$ iff for each  pair of indices $i\neq j$ where $1\leq i,j\leq s= |\maxDen_l(\bR ) |$ there exist $S_i, S_j\in \Den_l(R)$ such that $0\in S_iS_j$ (where $S_iS_j$ is the multiplicative  submonoid of $R$ generated by $S_i$ and $S_j$), $p_i(S_i) \subseteq \bQ_i^*$ and $p_j(S_j) \subseteq \bQ_j^*$.
 In this case,
 $\maxDen_l(R)=\{S_1', \ldots , S_s'\}$ where $S_i'$ is the largest element in the set $\{ S' \in \Den_l(R)\, | \, p_i(S')\subseteq \bQ_i^*\}$.
\end{theorem}

{\it Proof}. By \cite[Theorem 4.1]{Bav-Crit-S-Simp-lQuot}, $\maxDen_l(\bR ) = \{ T_1, \ldots , T_s\}$ where $T_i= (\bp_i\s)^{-1}(\bQ_i^*)$ for $i=1, \ldots , s$. 

$(\Rightarrow )$ It suffices to take $S_i=T_i$,  $i=1, \ldots , s$ since $0\in S_iS_j$ for all $i\neq j$, by Lemma \ref{a14Sep13}.(2).

$(\Leftarrow )$ Suppose that $S_1, \ldots , S_s$ are as in the theorem. For each $i=1, \ldots , s$, let $S_i'$ be the {\em largest} element with respect to  inclusion of the set $\{ S'\in \Den_l(R)\, | \, \pi (S')\subseteq T_i\}$. Since $S_i\subseteq  S_i'$ and $S_j\subseteq S_j'$ for each distinct pair $i\neq j$ and $0\in S_iS_j$, the elements $S_1',\ldots , S_s'$ are distinct  (Lemma \ref{a14Sep13}.(2)) and incomparable (i.e. $S_i'\not\subseteq S_j'$ for all $i\neq j$). In the proof of Theorem \ref{aA5Oct14}, the statement (vi), we have seen that $\maxDen_l(R)\subseteq \{ S_1',\ldots , S_s'\}$. Hence, $\maxDen_l(R)=\{ S_1',\ldots , S_s'\}$  since the sets
 $S_1',\ldots , S_s'$ are  incomparable.  $\Box $

$\noindent $

The next theorem gives sufficient conditions for finiteness of the set of $\maxDen_l(R)$.

\begin{theorem}\label{3Jul15}
Let $R$ be a ring and $I$ be an ideal of the ring $Q_l(R)$ such that $S+I\subseteq S$ and $S^{-1}I$ is an ideal of $S^{-1}Q_l(R)$ for all $S\in \maxDen_l(Q_l(R))$. Then $|\maxDen_l(R)|\leq |\maxDen_l(Q_l(R)/I)|$. In particular, if $|\maxDen_l(Q_l(R)/I)|<\infty$ then $|\maxDen_l(R)|<\infty$.
\end{theorem}

{\it Proof}. By Proposition \ref{A8Dec12}.(4), there is a bijection between the sets $\maxDen_l(R)$ and $\maxDen_l(Q_l(R))$. Let $Q_l:=Q_l(R)$ and $f:Q_l\ra \bQ_l:=Q_l/I$, $ q\mapsto \bq:=q+I$. 

(i) {\em For all} $S\in \maxDen_l(Q_l)$, $\bS:=f(S)\in \Den_l(\bQ_l)$: Since $S+I\subseteq S$, $S\cap I=\emptyset$ and $\bS\in \Ore_l(\bR )$. It remains to show that if $\bq\bs =0$ for some elements $\bq\in \bQ_l$ and $\bs \in \bS$ then $\bt\bq =0$ for some element $\bt \in \bS$. The equality $\bq\bs =0$ means that $i:=qs\in I$. Then $S^{-1}Q_l\ni \frac{q}{1}= is^{-1} = s_1^{-1}i_1$ for some elements $s_1\in S$ and $i_1\in I$ (since $S^{-1}R$ is an ideal of the ring $S^{-1}Q_l$). Hence, $s_2s_1q=s_2i\in I$ for some $s_2\in S$, and so $\bt\bq =0$ where $\bt =\bs_2\bs_1\in \bS$.

(ii)  $|\maxDen_l(R)|\leq |\maxDen_l(Q_l(R)/I)|$: Let $S\in \maxDen_l(Q_l)$. By (i), $\bS\subseteq T$ for some       $T=T(S)\in \maxDen_l(\bQ_l)$. Then $S$ is the largest element (w.r.t. $\subseteq $) of the set $\{ S'\in \Den_l(R)\, | \, f(S')\subseteq T\}$ (as $ 0\in SS_1$   for all distinct $S, S_1\in \maxDen_l(R)$, by Lemma \ref{a14Sep13}.(2)), and the statement (ii) follows.
$\Box $


\begin{corollary}\label{a3Jul15}
Let $R$ be a ring and $I$ be an ideal of $Q_l(R)$ such that $S^{-1}I$ is an ideal of the ring $S^{-1}Q_l(R)$ with $S^{-1}I\subseteq \rad (S^{-1}Q_l(R))$. Then $|\maxDen_l(R)|\leq |\maxDen_l(Q_l(R)/I)|$.
\end{corollary}

{\it Proof}. In view of Theorem \ref{3Jul15} and its proof, it suffices to show that $S+I\subseteq S$ for all $S\in \maxDen_l(Q_l)$ where $Q_l = Q_l(R)$. By the assumption, $S^{-1}I\subseteq \rad (S^{-1}Q_l)$. Hence, $1+S^{-1} I\subseteq (S^{-1}Q_l)^*$. Now, for all $s\in S$ and $i\in I$, $s+i=s(1+s^{-1}i)\in (S^{-1}Q_l)^*$. Let $\s_S:Q_l\ra S^{-1}Q_l$, $q\mapsto \frac{q}{1}$. By Theorem \ref{4Jul10}.(1), $S=\s_S^{-1} ((S^{-1}Q_l)^*)$, and so $s+i\in S$, as required.
$\Box $

$\noindent $

{\bf Proof of Theorem \ref{23Nov14} (description of the set $\maxDen_l(R)$ for a commutative ring $R$)}.
 For each minimal prime ideal $\gp$ of $R$, $S_\gp \in \maxDen (R)$ since the ring $R_\gp := S_\gp^{-1}R$ is a local ring such that its  maximal ideal $S^{-1}\gp$ (which is also the prime radical of $R_\gp$)  is a nil ideal, so every element of $R_\gp $ is either a unit or a nilpotent element.
 
 Conversely, let $S\in \maxDen (R)$. We have to show that $S= S_\gp$ for some $\gp \in \Min (R)$. In the ring $S^{-1}R$, every element is either a unit or a nilpotent element. Hence, the ring $S^{-1}R
 $ is a local ring $(S^{-1}R, \gm )$ where the maximal ideal $\gm$ is a nil ideal. Let $\s : R\ra S^{-1}R$, $r\mapsto \frac{r}{1}$. Then $S=\s^{-1} ((S^{-1}R)^*)= S_\gp$, by Theorem \ref{15Nov10}.(3).  $\Box $


\begin{proposition}\label{a22Nov14}
Let $R$ be a commutative ring, $I$ be a nil ideal of $R$ and $\pi_I:R\ra R/I$, $ r\mapsto \br = r+I$. Then \begin{enumerate}
\item The map $\Den(R)\ra  \Den (R/ I)$, $S\mapsto \pi_I(S)$, is a surjection that respects inclusions. Moreover, the map $\Den(R, I):=\{ S\in \Den (R)\, | \, S+I\subseteq S\} \ra  \Den (R/ I)$, $S\mapsto \pi_I(S)$, is a bijection with the inverse $T\mapsto \pi^{-1}_I(T)$.
\item The map $\maxDen(R)\ra  \maxDen (R/ I)$, $S\mapsto \pi_I(S)$, is a bijection with the inverse $T\mapsto \pi^{-1}_I(T)$.
\end{enumerate}
\end{proposition}

{\it Proof}. 1. Straightforward.

2. Statement 2 follows statement 1 since if $S\in \maxDen (R)$ then $S+I\subseteq S$. $\Box$

Let $R$ be a {\em commutative} ring such that its  quotient ring $Q= Q(R)$ is a {\em Noetherian} ring. We keep the notation of Theorem \ref{a5Oct14} and Theorem \ref{6Jan13}. For each $i=1, \ldots , s$, let $p_i:R\ra \bQ_i$ be as in (\ref{piRRQ}).
The next theorem shows that the second inequality of Theorem \ref{aA5Oct14} is an equality, it also provides another characterization of the set $\maxDen (R)$.

\begin{theorem}\label{b22Nov14}
Let $R$ be a commutative ring such that its quotient ring $Q$ is a Noetherian ring. Then $\maxDen (R) = \{ S_i\, | \, i=1, \ldots , s\}$ where $S_i := p_i^{-1} (\bQ_i^*)$. In particular, $|\maxDen(R)|=s=\maxDen (\bR )$.
\end{theorem}

{\it Proof}. Notice that $\maxDen (\bQ ) = \{ \bp_i^{-1} (\bQ_i^*)\, | \, i=1, \ldots , s\}$. By Proposition \ref{A8Dec12}.(4), the map
$\maxDen (\bQ) \ra \maxDen (\bR)$, $ S\mapsto \overline{\s}^{-1} (S)$, is a bijection. By Proposition \ref{a22Nov14}.(2), the map $\maxDen (\bR) \ra \maxDen (R)$, $ T\mapsto \pi^{-1} (T)$, is a bijection. Therefore, $\maxDen (R) = \{ S_i\, | \, i=1, \ldots , s\}$. $\Box $


\begin{proposition}\label{A9Nov14}
Let $R$ be a ring, $S\in \Den_l(R, \ga )$, $I$ be an ideal of $R$ such that $I\cap S=\emptyset$ and $S^{-1}I$ is an ideal of $S^{-1}R$, $\pi : R\ra R/ I$, $r\mapsto \br:= r+I$, and $\bS := \pi (S)$. Then
\begin{enumerate}
\item $\bS\in \Den_l(R/I, \overline{\gb} )$ for some ideal $\overline{\gb}$ of $R/I$ such that $\ga +I\subseteq \overline{\gb}$.
\item Let $\gb := \pi^{-1}(\overline{\gb})$. Then the map $\phi : S^{-1}R/ S^{-1}I\ra \bS^{-1} (R/I)$, $ s^{-1}r+S^{-1}I\mapsto \bs^{-1}\br$, is an epimorphism with $\ker (\phi ) = S^{-1}\gb / S^{-1}I$.
\end{enumerate}
\end{proposition}

{\it Proof}. 1. Since $I\cap S=\emptyset$, $\bS\in \Ore_l(R/I)$. It remains to show that if $\br\bs =0$ for some elements $\br\in R/I$ and $\bs \in \bS$ then $\bt \br =0$ for some $\bt\in \bS$.
 The equality $\br \bs=0$ means that $i:=rs\in I$. Hence, $S^{-1}R\ni \frac{r}{1}=is^{-1} = s_1^{-1} i_1$ for some elements $s_1\in S$ and $i_1\in I$ (since $S^{-1}I$ is an ideal of $S^{-1}R$). Therefore, $s_2s_1r=s_2i_1$ for some element $s_2\in S$, and so $\bt \br =0$ where $\bt = \bs_2\bs_1\in \bS$. Clearly, $\ga +I\subseteq \overline{\gb}$.

2. It is obvious that $\phi$ is an epimorphism. Since $\bS^{-1} (R/I) \simeq S^{-1}(R/I)$ and $I\subseteq \gb$, we have $\ker (\phi )=  S^{-1}\gb / S^{-1}I$.  $\Box $


\begin{proposition}\label{A8Nov14}
Let $R$ be a ring, $S\in \Den_l(R, 0)$, $Q':= S^{-1}R$ and $Q'^*$ be the group of units of $Q'$. Suppose that $S=R\cap Q'^*$ and one of the following statements holds:
\begin{enumerate}
\item Every left invertible element of $Q'$ is invertible.
\item Every right invertible element of $Q'$ is invertible.
\item The ring $Q'$ is a  left Noetherian ring.
\item The ring $Q'$ is a  right Noetherian ring.
\item The ring $Q'$ satisfies ACC on left annihilators.
\item The ring $Q'$ satisfies ACC on right annihilators.
\item The ring $Q'$ does not contain infinite direct sums of nonzero left ideals.
\item The ring $Q'$ does not contain infinite direct sums of nonzero right ideals.
\end{enumerate}
If $yz\in S$ for some elements $y,z\in R$ then $y,z\in S$.
\end{proposition}

{\it Proof}. If  $q:=yz\in S$ then $q\in Q'^*$. In particular, $y\cdot zq^{-1} = 1$ and $ q^{-1} y \cdot z=1$. If condition 1 (resp. 2) holds then $z$ (resp. $y$) is a unit of $Q'^*$. So, if any conditions 1 or 2 holds then $y,z\in Q'^*$, and so $y,z\in R\cap Q'^* = S$. Clearly, $(3\Rightarrow 5)$ and $(4\Rightarrow 6)$. So, it suffices to consider the case when one of the conditions 5-8 holds. Let $x:= zq^{-1}$. Then $yx=1$ in $Q'$. The idempotents $e_i:= x^iy^i - x^{i+1}y^{i+1}= x^i e_0y^i$, $i=0,1,\ldots $, are orthogonal idempotents. Since one of  conditions 5-8 holds we must have $e_i=0$ for some $i$ then $0=y^ie_ix^i = y^ix^i e_0y^ix^i = e_0= 1-xy$, i.e. $1=xy = yx$. Hence, $x,y\in Q'^*$ and $x,y\in R\cap Q'^* = S$.  $\Box $


\begin{corollary}\label{aA8Nov14}
\begin{enumerate}
\item Let $R$ be a ring such that its largest left quotient ring $Q':= Q_l(R)$ satisfies one of conditions 1-8 of Proposition \ref{A8Nov14}. If $yz\in S_l(R)$ for some elements $y,z\in R$ then $y,z\in S_l(R)$.
\item Let $R$ be a ring such that its  left quotient ring $Q':= Q(R)$ satisfies one of conditions 1-8 of Proposition \ref{A8Nov14}. If $yz\in \CC_R$ for some elements $y,z\in R$ then $y,z\in \CC_R$.
\end{enumerate}
\end{corollary}

{\it Proof}. 1. Statement 1 follows from Proposition \ref{A8Nov14} since $S_l(R) = R\cap Q_l(R)^*$ (Theorem \ref{4Jul10}.(1)).

2. Statement 2 is a particular case of statement 1.  $\Box $

$\noindent $



\small{

Department of Pure Mathematics

University of Sheffield

Hicks Building

Sheffield S3 7RH

UK

email: v.bavula@sheffield.ac.uk}

\end{document}